\documentclass[namedreferences,letterpaper,12pt]{kluwer}
\usepackage[left=1in, right=1in, top=1in, bottom=1.25in]{geometry}
\usepackage{pgf,pgfarrows,pgfnodes} 
\usepackage{times}
\usepackage{fancyhdr}
\begin{document}
\begin{opening} 
\title{The Uses of Argument in Mathematics}
\def\authorsize{\large \raggedright}
\def\afterauthorskip{\aftertitleskip}
\author{ANDREW \surname{ABERDEIN}} 
\institute{Humanities and Communication,\\
Florida Institute of Technology,\\
150 West University Blvd, \\
Melbourne, Florida 32901-6975, U.S.A.\\
aberdein@fit.edu}
\renewcommand{\abstractname}{ABSTRACT}
\renewcommand{\keywordsname}{KEY WORDS}
\def\abstractnamefont{\rm}
\def\keynamefont{\rm}
\begin{abstract} 
Stephen Toulmin once observed that `it has never been customary for philosophers to pay much attention to the rhetoric of mathematical debate' \protect\cite[p.~89]{Toulmin+}.  Might the application of Toulmin's layout of arguments to mathematics remedy this oversight?

Toulmin's critics fault the layout as requiring so much abstraction as to permit incompatible reconstructions.  Mathematical proofs may indeed be represented by fundamentally distinct layouts.  However, cases of genuine conflict characteristically reflect an underlying disagreement about the nature of the proof in question. 
 \end{abstract} 
\end{opening} 

\setlength{\parindent}{.5in}
\pagestyle{fancy}
\chead{\sc Andrew Aberdein, The Uses of Argument in Mathematics} 
\renewcommand{\headrulewidth}{0pt} 
\makeatletter 
\renewcommand{\thesubsection}{\@arabic\c@subsection} 
\renewcommand{\@copyrightfoot}{}
\renewcommand{\@copyrighthead}{}
\renewcommand{\idline}{}
\makeatother 

\noindent Modern formal logic's earliest successes were in the analysis of mathematical argument.  However, formal logic has long been recognized as useful elsewhere, and is now the dominant means of logical analysis in all fields of discourse.  Informal logicians claim that this mathematical success has been bought at the cost of broader application.  They have proposed methods of argument analysis complementary to that of formal logic, providing for the pragmatic treatment of features of argumentation which cannot be reduced to logical form.  Characteristically, informal logicians have conceded mathematics to formal logic, while stressing the superiority of their own systems as analyses of argumentation in natural language.  My contention is that this is an unnecessary concession: both systems have lessons for mathematics, just as they do for natural language argumentation.

In this paper I exhibit some aspects of mathematical argumentation which can best be captured by informal logic.  Specifically I investigate the applicability of Toulmin's layout of arguments to mathematics.  I demonstrate how the layout may be used to represent the structure of both `regular' and `critical' arguments in mathematics.  
Mathematical proof is typically regular argumentation, and it is this which I address most closely.

\subsection{Toulmin's Layout}

Toulmin's \textit{The Uses of Argument} is perhaps the single most influential work in modern argumentation theory \cite{Toulmin58}.  Especially widely cited is the general account it offers of the structure of arguments.  Toulmin begins with the thought that an argument is a claim ($C$) derived from data ($D$) in accordance with a warrant ($W$).  While this is superficially similar to the treatment of arguments in deductive logic, greater generality is achieved through additional components of the layout.  The argument may have a modal qualifier ($Q$), such as `necessarily' or `presumably', which explicates the force of the warrant.  If the warrant does not provide necessity, its conditions of exception or rebuttal ($R$) may be noted.  We may also keep track of the backing ($B$) which supports the warrant.  The overall layout is often set out graphically, as in Figure~\ref{fig:Toulmin}.

\begin{figure}[h]
\label{fig:Toulmin}
\centering
\begin{pgfpicture}{0in}{0in}{2in}{1.7in}
\pgfsetlinewidth{.5pt} 
\pgfnodebox{NodeQ1}[stroke]
   {\pgfxy(2.6,1.5)}
   {Q}{2pt}{2pt} 
\pgfnodebox{NodeD1}[stroke]
   {\pgfrelative{\pgfxy(-1.5,0)}{\pgfnodeborder{NodeQ1}{180}{0pt}}}
   {D}{2pt}{2pt}  
\pgfnodebox{NodeC1}[stroke]
   {\pgfrelative{\pgfxy(1.5,0)}{\pgfnodeborder{NodeQ1}{0}{0pt}}}
   {C}{2pt}{2pt}  
\pgfnodebox{NodeR1}[stroke]
   {\pgfrelative{\pgfxy(0,-1)}{\pgfnodeborder{NodeQ1}{270}{0pt}}}
   {R}{2pt}{2pt}  
\pgfnodebox{Invisible1}[virtual]
   {\pgfrelative{\pgfxy(-.6,0)}{\pgfnodeborder{NodeQ1}{180}{0pt}}}{}{0pt}{0pt} 
\pgfnodebox{NodeW1}[stroke]
   {\pgfrelative{\pgfxy(0,1)}{\pgfnodecenter{Invisible1}}}
   {W}{2pt}{2pt}
\pgfnodebox{NodeB1}[stroke]
   {\pgfrelative{\pgfxy(0,1)}{\pgfnodeborder{NodeW1}{90}{0pt}}}
   {B}{2pt}{2pt}

\pgfsetlinewidth{1pt} 
\pgfnodeconnline{NodeD1}{NodeQ1} 
\pgfsetendarrow{\pgfarrowsingle} 
\pgfnodesetsepstart{0pt}\pgfnodesetsepend{2pt} 
\pgfnodeconnline{NodeQ1}{NodeC1} 
\pgfnodeconnline{NodeR1}{NodeQ1} 
\pgfnodeconnline{NodeW1}{Invisible1} 
\pgfnodeconnline{NodeB1}{NodeW1} 
\end{pgfpicture}
\caption{Toulmin's layout \protect\cite[p.~104]{Toulmin58}.}
\end{figure}

While philosophers attacked \textit{The Uses of Argument} as `Toulmin's anti-logic book', rhetoricians, communication theorists and subsequently computer scientists have found its account of argumentation deeply insightful \cite{Shapin}.  Although it traduces Toulmin to characterize him as `anti-logic', his motivation was indeed to critique formal logic, specifically its claim to represent all the significant content of \textit{non}-mathematical argument.  Mathematics, as the discourse for which formal logic is best suited, was exempt from this critique.
Nevertheless, Toulmin's layout is intended to encompass all forms of argument, mathematics included, as he makes clear in a later work \protect\cite[p.~89]{Toulmin+}.  In the remainder of this paper we will see how well it accomplishes this task, and whether it has advantages over formal logic alone in this discourse too.

One significant problem which critics of Toulmin have raised is that the degree of abstraction necessary to use the layout at all can make different, incompatible, reconstructions possible \cite{Willard}.  This would seem to imply that use of the layout must risk distorting the original argument to an unacceptable degree.  I shall explore whether this criticism undermines the application of Toulmin's layout to mathematics below, but first I must address an important distinction.

\subsection{Regular and Critical Arguments}
\label{sec:Reg}
In discussing scientific arguments, Toulmin distinguishes between arguments that are conducted within, or as applications of, a scientific theory and arguments which challenge a prevailing theory or seek to motivate an alternative \cite[p.~247]{Toulmin+}.  The former he terms `regular arguments', the latter `critical arguments'.  This distinction echoes the familiar one between normal and revolutionary science promoted by Thomas Kuhn, although Toulmin rejects Kuhn's hard and fast distinction between these two modes \cite{Toulmin70}.  Hence critical arguments are always available to scientists, not just in revolutionary phases, although they may not always have cause to use them.  Moreover, critical arguments are more broadly conceived than scientific revolutions, since they do not necessarily require the overthrow of the prior theory in order to succeed.  
As Toulmin observes, `[c]omplete immutability in our rational procedures and standards of judgment is not to be found even in \dots pure mathematics. The standards of rigour relied on in the judging of mathematical arguments have had their own history' \cite[p.~133]{Toulmin+}.
Even commentators who deny that there can be revolutions in mathematics on the grounds that it is a purely cumulative discipline concede that there can be revolutions in mathematical rigour.\footnote{For example, \cite[p.~19]{Crowe}.  Positions for and against mathematical revolutions are explored at length in \cite{Gillies}.}
Thus even the most conservative history of mathematics must consider mathematical critical arguments.  Furthermore, when we apply Toulmin's layout to regular mathematical argumentation, we shall see how it may be used to keep track of the changing standards of rigour to which critical argumentation can give rise.

\begin{figure}[h]
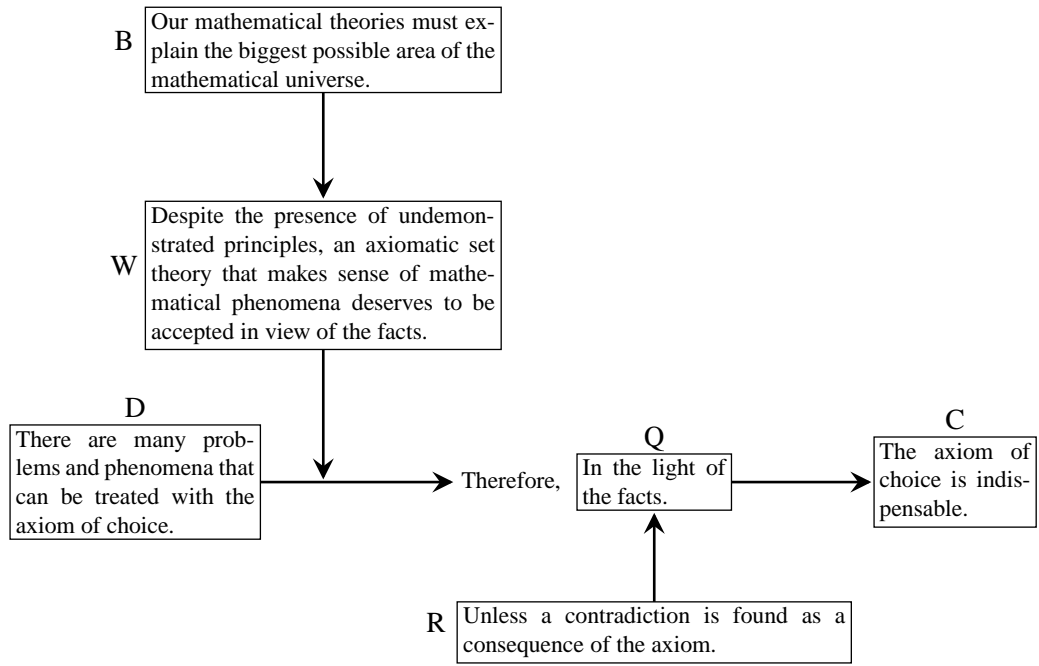
\label{fig:Alcolea}
\begin{center}
\begin{pgfpicture}{-0.3in}{-3in}{5in}{.7in}
\pgfsetlinewidth{.5pt} 
\pgfnodebox{NodeD}[stroke]
   {\pgfxy(1,-5)}
   {\parbox{1.25in}{\footnotesize There are many problems and phenomena that can be treated with the axiom of choice.}}{2pt}{2pt} 
\pgfnodebox{Therefore}[virtual]
   {\pgfrelative{\pgfxy(5,0)}{\pgfnodecenter{NodeD}}}
   {\footnotesize Therefore,}{2pt}{2pt} 
\pgfnodebox{Invisible}[virtual]
   {\pgfrelative{\pgfxy(2.5,0)}{\pgfnodecenter{NodeD}}}{}{0pt}{0pt} 
\pgfnodebox{NodeW}[stroke]
   {\pgfrelative{\pgfxy(0,2.75)}{\pgfnodecenter{Invisible}}}
   {\parbox{1.8in}{\footnotesize Despite the presence of undemonstrated principles, an axiomatic set theory that makes sense of mathematical phenomena deserves to be accepted in view of the facts.}}{2pt}{2pt}
\pgfnodebox{NodeB}[stroke]
   {\pgfrelative{\pgfxy(0,3)}{\pgfnodecenter{NodeW}}}
   {\parbox{1.8in}{\footnotesize Our mathematical theories must explain the biggest possible area of the mathematical universe.}}{2pt}{2pt}
\pgfnodebox{NodeQ}[stroke]
   {\pgfrelative{\pgfxy(1.9,0)}{\pgfnodecenter{Therefore}}}
   {\parbox{0.75in}{\footnotesize In the light of the facts.}}{2pt}{2pt}  
\pgfnodebox{NodeC}[stroke]
   {\pgfrelative{\pgfxy(4,0)}{\pgfnodecenter{NodeQ}}}
   {\parbox{.79in}{\footnotesize The axiom of choice is indispensable.}}{2pt}{2pt}
\pgfnodebox{NodeR}[stroke]
   {\pgfrelative{\pgfxy(0,-2)}{\pgfnodecenter{NodeQ}}}
   {\parbox{2in}{\footnotesize Unless a contradiction is found as a consequence of the axiom.}}{2pt}{2pt}

\pgfputat{\pgfnodeborder{NodeD}{90}{3pt}}{\pgfbox[center,base]{D}}
\pgfputat{\pgfnodeborder{NodeQ}{90}{3pt}}{\pgfbox[center,base]{Q}}
\pgfputat{\pgfnodeborder{NodeC}{90}{3pt}}{\pgfbox[center,base]{C}}
\pgfputat{\pgfnodeborder{NodeR}{180}{8pt}}{\pgfbox[center, center]{R}}
\pgfputat{\pgfnodeborder{NodeW}{180}{8pt}}{\pgfbox[center, center]{W}}
\pgfputat{\pgfnodeborder{NodeB}{180}{8pt}}{\pgfbox[center, center]{B}}

\pgfsetlinewidth{1pt} 
\pgfsetendarrow{\pgfarrowsingle} 
\pgfnodesetsepstart{0pt}\pgfnodesetsepend{2pt} 
\pgfnodeconnline{NodeD}{Therefore} 
\pgfnodeconnline{NodeQ}{NodeC} 
\pgfnodeconnline{NodeR}{NodeQ} 
\pgfnodeconnline{NodeB}{NodeW} 
\pgfnodeconnline{NodeW}{Invisible} 
\end{pgfpicture} 

\end{center}
\caption{Alcolea's analysis of Zermelo's argument for the adoption of the axiom of choice \protect\cite[p.~143]{Alcolea}.}
\end{figure}

A good example of a mathematical critical argument is that offered by Ernst Zermelo and others for admitting the Axiom of Choice as one of the axioms of set theory.  This is discussed by Jes\'{u}s Alcolea Banegas, one of the few people to specifically address the application of informal logic to mathematics \cite{Alcolea}.\footnote{I am grateful to Miguel Gimenez of the University of Edinburgh for translating this paper from the original Catalan.}  I have reproduced his reconstruction of the layout of this argument as Figure~\ref{fig:Alcolea}.
However, mathematical critical arguments are arguments \textit{about} mathematics, not arguments \textit{in} mathematics,  For this reason, mathematical critical arguments are much like the critical arguments of any other discipline: there is nothing specifically mathematical about the warrant invoked, or any other aspect of their structure.  Indeed, there cannot be, since they are appealing to the extra-mathematical in order to settle a dispute which cannot be resolved by purely mathematical means.

The characteristic content of mathematics is mathematical proof.  Proofs are regular arguments---they may inspire criticism of underlying principles, but the criticism must take place outside the proof itself.  In the next section we will examine how closely Toulmin's layout models the work of mathematical proof.

\subsection{Proofs}

Toulmin's own example applying his layout to mathematics is Theaetetus's proof that there are exactly five platonic solids.  This result is recorded by Euclid as the final proposition of Book XIII of his \textit{Elements}, where the proof reads:

\begin{quotation}
I say next that \textit{no other figure, besides the said five figures, can be constructed which is contained by equilateral and equiangular figures equal to one another}.

For a solid angle cannot be constructed with two triangles, or indeed planes.

With three triangles the angle of the pyramid is constructed, with four the angle of the octahedron, and with five the angle of the icosahedron; but a solid angle cannot be formed by six equilateral and equiangular triangles placed together at one point, for, the angle of the equilateral triangle being two thirds of a right angle, the six will be equal to four right angles: which is impossible, for any solid angle is contained by angles less than four right angles [XI.21].

For the same reason, neither can a solid angle be constructed by more than six plane angles.

By three squares the angle of the cube is contained, but by four it is impossible for a solid angle to be contained, for they will again be four right angles.

By three equilateral and equiangular pentagons the angle of the dodecahedron is contained; but by four such it is impossible for any solid angle to be contained, for, the angle of the equilateral pentagon being a right angle and a fifth, the four angles will be greater than the four right angles: which is impossible.

Neither again will a solid angle be contained by other polygonal figures by reason of the same absurdity.

\raggedleft{\sc q.e.d.} \cite[Vol.~3, pp.~507 f.]{Heath}
\end{quotation}

\begin{figure}[h]
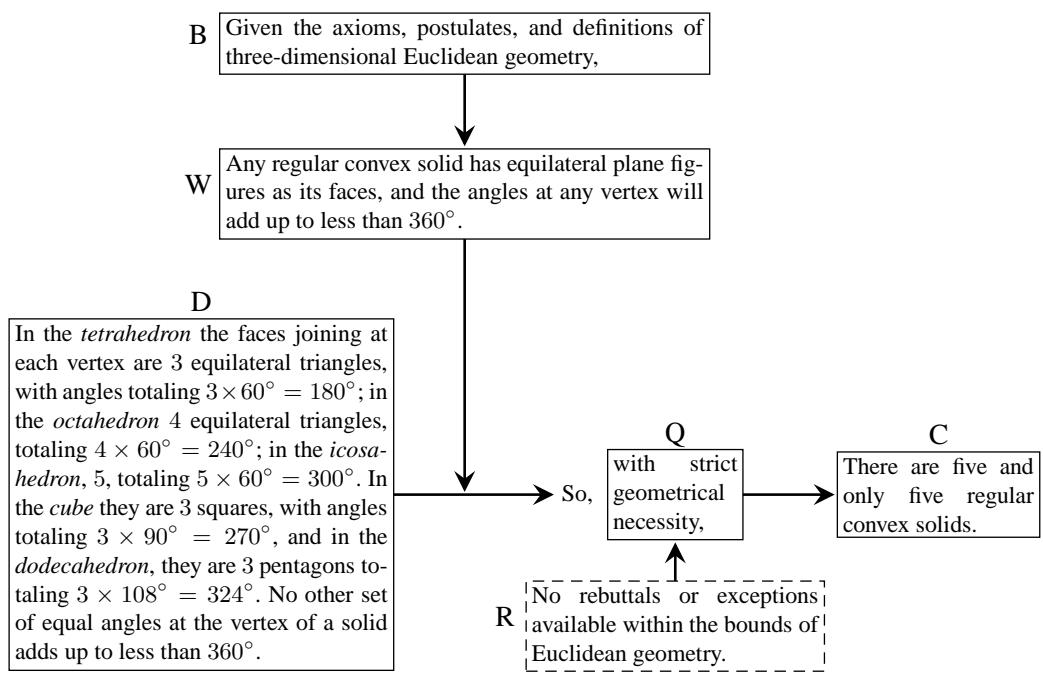
 \label{fig:Theaetetus}
\begin{center}
\begin{pgfpicture}{-0.3in}{-3.2in}{5in}{.5in}
\pgfsetlinewidth{.5pt} 
\pgfnodebox{NodeD}[stroke]
   {\pgfxy(1.8,-5.5)}
   {\parbox{1.95in}{\footnotesize In the \textit{tetrahedron} the faces joining at each vertex are $3$ equilateral triangles, with angles totaling $3 \times 60^{\circ} = 180^{\circ}$; in the \textit{octahedron} $4$ equilateral triangles, totaling $4 \times 60^{\circ} = 240^{\circ}$; in the \textit{icosahedron}, $5$, totaling $5 \times 60^{\circ} = 300^{\circ}$.  In the \textit{cube} they are $3$ squares, with angles totaling $3 \times 90^{\circ} = 270^{\circ}$, and in the \textit{dodecahedron}, they are $3$ pentagons totaling $3 \times 108^{\circ} = 324^{\circ}$.  No other set of equal angles at the vertex of a solid adds up to less than $360^{\circ}$.}}{2pt}{2pt} 
\pgfnodebox{Therefore}[virtual]
   {\pgfrelative{\pgfxy(5,0)}{\pgfnodecenter{NodeD}}}
   {\footnotesize So,}{2pt}{2pt} 
\pgfnodebox{Invisible}[virtual]
   {\pgfrelative{\pgfxy(3.5,0)}{\pgfnodecenter{NodeD}}}{}{0pt}{0pt} 
\pgfnodebox{NodeW}[stroke]
   {\pgfrelative{\pgfxy(0,4)}{\pgfnodecenter{Invisible}}}
   {\parbox{2.5in}{\footnotesize Any regular convex solid has equilateral plane figures as its faces, and the angles at any vertex will add up to less than $360^{\circ}$.}}{2pt}{2pt}
\pgfnodebox{NodeB}[stroke]
   {\pgfrelative{\pgfxy(0,2)}{\pgfnodecenter{NodeW}}}
   {\parbox{2.5in}{\footnotesize Given the axioms, postulates, and definitions of three-dimensional Euclidean geometry,}}{2pt}{2pt}
\pgfnodebox{NodeQ}[stroke]
   {\pgfrelative{\pgfxy(1.3,0)}{\pgfnodecenter{Therefore}}}
   {\parbox{0.65in}{\footnotesize with strict geometrical necessity,}}{2pt}{2pt}  
\pgfnodebox{NodeC}[stroke]
   {\pgfrelative{\pgfxy(3.5,0)}{\pgfnodecenter{NodeQ}}}
   {\parbox{1in}{\footnotesize There are five and only five regular convex solids.}}{2pt}{2pt}
\pgfsetdash{{0.06in}{0.04in}}{0in}
\pgfnodebox{NodeR}[stroke]
   {\pgfrelative{\pgfxy(0,-1.75)}{\pgfnodecenter{NodeQ}}}
   {\parbox{1.5in}{\footnotesize No rebuttals or exceptions available within the bounds of Euclidean geometry.}}{2pt}{2pt}
\pgfsetdash{{1in}{0in}}{0in}

\pgfputat{\pgfnodeborder{NodeD}{90}{3pt}}{\pgfbox[center,base]{D}}
\pgfputat{\pgfnodeborder{NodeQ}{90}{3pt}}{\pgfbox[center,base]{Q}}
\pgfputat{\pgfnodeborder{NodeC}{90}{3pt}}{\pgfbox[center,base]{C}}
\pgfputat{\pgfnodeborder{NodeR}{180}{8pt}}{\pgfbox[center,base]{R}}
\pgfputat{\pgfnodeborder{NodeW}{180}{8pt}}{\pgfbox[center,base]{W}}
\pgfputat{\pgfnodeborder{NodeB}{180}{8pt}}{\pgfbox[center,base]{B}}

\pgfsetlinewidth{1pt} 
\pgfsetendarrow{\pgfarrowsingle} 
\pgfnodesetsepstart{0pt}\pgfnodesetsepend{2pt} 
\pgfnodeconnline{NodeD}{Therefore} 
\pgfnodeconnline{NodeQ}{NodeC} 
\pgfnodeconnline{NodeR}{NodeQ} 
\pgfnodeconnline{NodeB}{NodeW} 
\pgfnodeconnline{NodeW}{Invisible} 
\end{pgfpicture} 

\end{center}
\caption{Toulmin's analysis of Theaetetus's proof that the platonic solids are exactly five in number \protect\cite[Fig.~7.4, p.~89]{Toulmin+}.}
\end{figure}

Toulmin's layout of this proof, slightly adapted for consistency of notation, is reproduced as Figure~\ref{fig:Theaetetus}.  We can see that it is an essentially faithful reproduction of Euclid's argument, but it will be profitable to compare the two more closely.  Firstly, there is a harmless simplificiation, but slight loss of rigour, in Toulmin's omission of the statement that a solid angle cannot be constructed from two planes (a consequence of Euclid's Def.\ 11, Book XI).  Most of the remainder of Euclid's proof becomes Toulmin's data, $D$, the claim, $C$, is essentially the same in both presentations, and Toulmin's warrant, $W$, combines the definition of a regular convex polyhedron with the essential prior result on which the proof depends, Euclid's Prop.\ 21, Book XI.  Other parts of the layout are more novel: $Q$, a characterization of the rigour of the proof, makes explicit something implicit in Euclid, and $R$ is not there at all, since on Toulmin's account the proof cannot be rebutted.  All of these identifications seem reasonable, and there does not seem to be any scope for pernicious ambiguity.  However, this is an elementary proof---problems may yet arise in more complicated cases.

One obvious way in which Theaetetus's proof is untypical of mathematical proofs is that it has only one step.  Most proofs, indeed most proofs in Euclid, constitute a sequence of steps, from the initial premisses, through one intermediate result after another to the eventual conclusion.  How might Toulmin's layout be applied to these multi-step proofs?  The most thorough method would be to diagram each step of the proof separately, perhaps linking them together in a sort of sorites.  If the fine detail of the proof is of particular interest this would be the best approach, but typically we would prefer something more coarse-grained.  A useful maxim that applies to any modeling process is that a model should be as detailed as we need it to be, \textit{but no more so}.  In most multi-step mathematical proofs the qualifier (and rebuttal) will be the same at every step.  This allows us to merge the layouts for the separate steps into one, by combining the given components of the data, the warrants (and if necessary the backing) of each step to produce a single layout for the whole proof.\footnote{This is a simplification, but more technical details are not required to establish this point.}  Where different steps have different qualifiers, the qualifier for the whole proof would represent the degree of certainty of the least certain step.  

\begin{figure}[h]
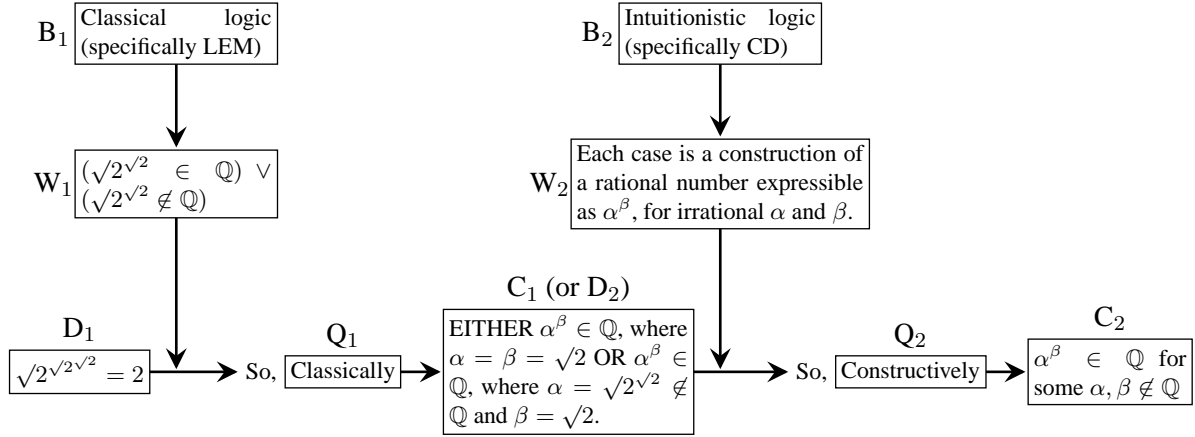

\begin{center}
\begin{pgfpicture}{.3in}{-2.5in}{5in}{0in}
\pgfsetlinewidth{.5pt} 
\pgfnodebox{NodeD}[stroke]
   {\pgfxy(-0.2,-5.3)}
   {{\footnotesize $\surd 2^{\surd 2^{\surd 2}}= 2$}}{2pt}{2pt} 
\pgfnodebox{Invisible}[virtual]
   {\pgfnodeborder{NodeD}{0}{10pt}}{}{0pt}{0pt} 
\pgfnodebox{NodeW}[stroke]
   {\pgfrelative{\pgfxy(0,2.5)}{\pgfnodecenter{Invisible}}}
   {\parbox{1in}{\footnotesize $(\surd 2^{\surd 2} \in \mathbb{Q}) \vee (\surd 2^{\surd 2} \not\in \mathbb{Q})$}}{2pt}{2pt}
\pgfnodebox{NodeB}[stroke]
   {\pgfrelative{\pgfxy(0,2)}{\pgfnodecenter{NodeW}}}
   {\parbox{1in}{\footnotesize Classical logic (specifically LEM)}}{2pt}{2pt}
\pgfnodebox{NodeC}[stroke]
   {\pgfrelative{\pgfxy(6.5,0)}{\pgfnodecenter{NodeD}}}
   {\parbox{1.25in}{\footnotesize EITHER $\alpha ^{\beta} \in  \mathbb{Q}$, where $\alpha = \beta = \surd 2$ OR $\alpha ^{\beta} \in  \mathbb{Q}$, where $\alpha = \surd 2^{\surd 2}  \not\in  \mathbb{Q}$ and $\beta = \surd 2$.}}{2pt}{2pt}
\pgfnodebox{NodeQ}[stroke]
   {\pgfrelative{\pgfxy(-3,0)}{\pgfnodecenter{NodeC}}}
   {{\footnotesize Classically}}{2pt}{2pt}  
\pgfnodebox{Therefore}[virtual]
   {\pgfnodeborder{NodeQ}{180}{9pt}}
   {\footnotesize So,}{2pt}{2pt} 
\pgfnodebox{Invisible2}[virtual]
   {\pgfnodeborder{NodeC}{0}{10pt}}{}{0pt}{0pt} 
\pgfnodebox{NodeW2}[stroke]
   {\pgfrelative{\pgfxy(0,2.5)}{\pgfnodecenter{Invisible2}}}
   {\parbox{1.5in}{\footnotesize 
   Each case is a construction of a rational number expressible as $\alpha ^{\beta}$, for irrational $\alpha$ and $\beta$.}}{2pt}{2pt}
\pgfnodebox{NodeB2}[stroke]
   {\pgfrelative{\pgfxy(0,2)}{\pgfnodecenter{NodeW2}}}
   {\parbox{1in}{\footnotesize Intuitionistic logic (specifically CD)}}{2pt}{2pt}
\pgfnodebox{NodeC2}[stroke]
   {\pgfrelative{\pgfxy(7.2,0)}{\pgfnodecenter{NodeC}}}
   {\parbox{.8in}{\footnotesize $\alpha ^{\beta}\in\mathbb{Q}$ for some $\alpha, \beta \not\in  \mathbb{Q}$}}{2pt}{2pt}
\pgfnodebox{NodeQ2}[stroke]
   {\pgfrelative{\pgfxy(-2.65,0)}{\pgfnodecenter{NodeC2}}}
   {{\footnotesize Constructively}}{2pt}{2pt}  
\pgfnodebox{Therefore2}[virtual]
   {\pgfnodeborder{NodeQ2}{180}{9pt}}
   {\footnotesize So,}{2pt}{2pt} 

\pgfputat{\pgfnodeborder{NodeD}{90}{5pt}}{\pgfbox[center,base]{D$_{1}$}}
\pgfputat{\pgfnodeborder{NodeW}{180}{8pt}}{\pgfbox[center,center]{W$_{1}$}}
\pgfputat{\pgfnodeborder{NodeB}{180}{8pt}}{\pgfbox[center,center]{B$_{1}$}}
\pgfputat{\pgfnodeborder{NodeQ}{90}{5pt}}{\pgfbox[center,base]{Q$_{1}$}}
\pgfputat{\pgfnodeborder{NodeC}{90}{5pt}}{\pgfbox[center,base]{C$_{1}$ (or D$_{2}$)}}
\pgfputat{\pgfnodeborder{NodeW2}{180}{8pt}}{\pgfbox[center,center]{W$_{2}$}}
\pgfputat{\pgfnodeborder{NodeB2}{180}{8pt}}{\pgfbox[center,center]{B$_{2}$}}
\pgfputat{\pgfnodeborder{NodeQ2}{90}{5pt}}{\pgfbox[center,base]{Q$_{2}$}}
\pgfputat{\pgfnodeborder{NodeC2}{90}{5pt}}{\pgfbox[center,base]{C$_{2}$}}

\pgfsetlinewidth{1pt} 
\pgfsetendarrow{\pgfarrowsingle} 
\pgfnodesetsepstart{0pt}\pgfnodesetsepend{2pt} 
\pgfnodeconnline{NodeD}{Therefore} 
\pgfnodeconnline{NodeQ}{NodeC} 
\pgfnodeconnline{NodeB}{NodeW} 
\pgfnodeconnline{NodeW}{Invisible} 
\pgfnodeconnline{NodeC}{Therefore2} 
\pgfnodeconnline{NodeQ2}{NodeC2} 
\pgfnodeconnline{NodeB2}{NodeW2} 
\pgfnodeconnline{NodeW2}{Invisible2} 
\end{pgfpicture} 

\end{center}
\caption{Classical proof that there are irrational numbers $\alpha$ and $\beta$ such that $\alpha ^{\beta}$ is rational.}
\label{fig:root2}\end{figure}

\begin{figure}[H]
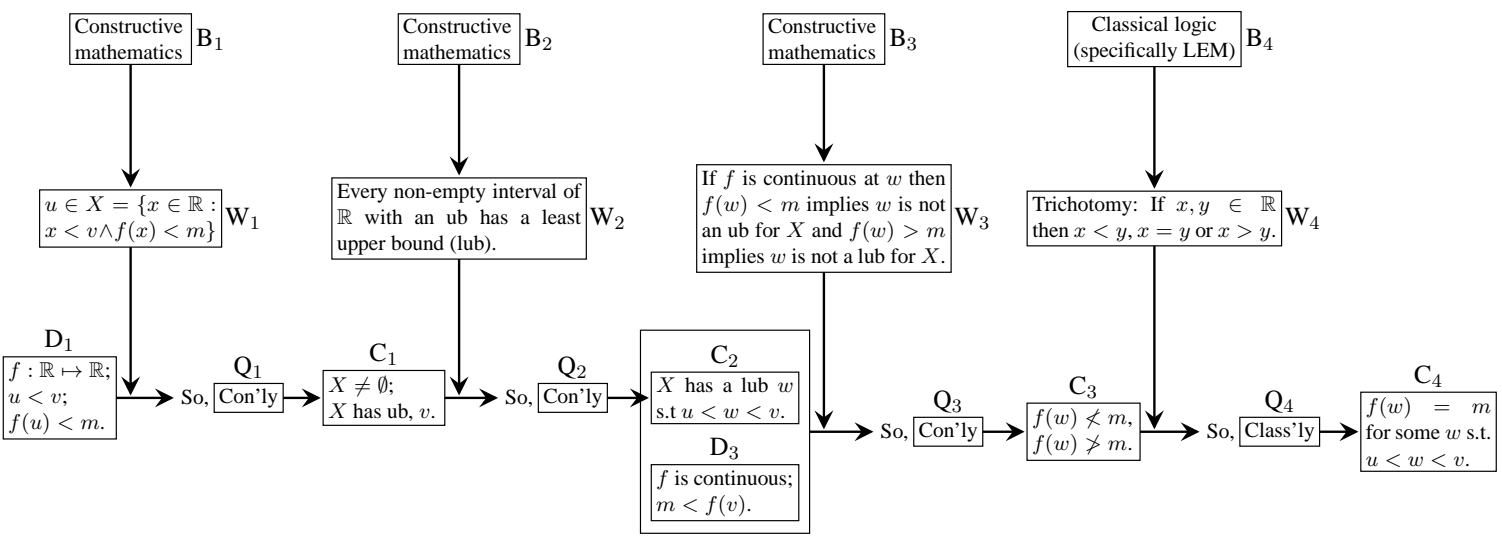
 
\begin{center}
\begin{pgfpicture}{.3in}{-8in}{5in}{0in}
\begin{pgfrotateby}{\pgfdegree{90}} 
\begin{pgfmagnify}{.9}{.9} 
\pgfsetlinewidth{.5pt} 
\pgfnodebox{NodeC3}[stroke]
   {\pgfxy(-11.75,-9.25)}
   {\parbox{.8in}{\footnotesize $X$ has a lub $w$ s.t $u<w<v$. 
   }}{2pt}{2pt} 
\pgfnodebox{NodeQ3}[stroke]
   {\pgfrelative{\pgfxy(-2.25,0)}{\pgfnodecenter{NodeC3}}}
   {{\footnotesize Con'ly}}{2pt}{2pt}  
\pgfnodebox{Therefore3}[virtual]
   {\pgfnodeborder{NodeQ3}{180}{8pt}}
   {\footnotesize So,}{2pt}{2pt} 
\pgfnodebox{NodeC1}[stroke]
   {\pgfrelative{\pgfxy(-2,0)}{\pgfnodecenter{Therefore3}}}
   {\parbox{.65in}
   {\footnotesize $X \neq \emptyset$;\\
    $X$ has ub, $v$.}}{2pt}{2pt}
\pgfnodebox{NodeQ1}[stroke]
   {\pgfrelative{\pgfxy(-2,0)}{\pgfnodecenter{NodeC1}}}
   {{\footnotesize Con'ly}}{2pt}{2pt}  
\pgfnodebox{Therefore1}[virtual]
   {\pgfnodeborder{NodeQ1}{180}{8pt}}
   {\footnotesize So,}{2pt}{2pt} 
\pgfnodebox{NodeD1}[stroke]
   {\pgfrelative{\pgfxy(-2,0)}{\pgfnodecenter{Therefore1}}}
   {\parbox{.6in}
   {\footnotesize 
   $f: \mathbb{R} \mapsto \mathbb{R}$;\\ 
   $u < v$;\\
   $f(u) < m$.}}{2pt}{2pt}  
\pgfnodebox{Invisible3}[virtual]
   {\pgfnodeborder{NodeC1}{0}{6pt}}{}{0pt}{0pt} 
\pgfnodebox{NodeW3}[stroke]
   {\pgfrelative{\pgfxy(0,2.65)}{\pgfnodecenter{Invisible3}}}
   {\parbox{1.42in}
   {\footnotesize 
   Every non-empty interval of $\mathbb{R}$ with an ub has a least upper bound (lub).}}{2pt}{2pt}
\pgfnodebox{NodeB3}[stroke]
   {\pgfrelative{\pgfxy(0,2.65)}{\pgfnodecenter{NodeW3}}}
   {\parbox{.65in}
   {\footnotesize Constructive mathematics}}{2pt}{2pt}
\pgfnodebox{Invisible1}[virtual]
   {\pgfnodeborder{NodeD1}{0}{6pt}}{}{0pt}{0pt} 
\pgfnodebox{NodeW1}[stroke]
   {\pgfrelative{\pgfxy(0,2.65)}{\pgfnodecenter{Invisible1}}}
   {\parbox{1in}
   {\footnotesize $u \in 
   X = \{x \in \mathbb{R}:
   x<v \wedge f(x)<m\}$}}{2pt}{2pt}
\pgfnodebox{NodeB1}[stroke]
   {\pgfrelative{\pgfxy(0,2.65)}{\pgfnodecenter{NodeW1}}}
   {\parbox{.65in}
   {\footnotesize Constructive mathematics}}{2pt}{2pt}

\pgfnodebox{NodeD3}[stroke]
   {\pgfrelative{\pgfxy(0,-1.4)}{\pgfnodecenter{NodeC3}}}
   {\parbox{.8in}
   {\footnotesize    $f$ is continuous;\\
    $m < f(v)$.}}{2pt}{2pt}
\pgfnoderect{NodeCD3}[stroke]{\pgfrelative{\pgfxy(0,-0.5)}{\pgfnodecenter{NodeC3}}}
{\pgfxy(2.5,3)} 
\pgfnodebox{Invisible4}[virtual]
   {\pgfnodeborder{NodeCD3}{0}{6pt}}{}{0pt}{0pt} 
\pgfnodebox{NodeW4}[stroke]
   {\pgfrelative{\pgfxy(0,3.15)}{\pgfnodecenter{Invisible4}}}
   {\parbox{1.42in}{\footnotesize 
   If $f$ is continuous at $w$ then $f(w) < m$ implies $w$ is not an ub for $X$ and $f(w) > m$ implies $w$ is not a lub for $X$.
}}{2pt}{2pt}
\pgfnodebox{NodeB4}[stroke]
   {\pgfrelative{\pgfxy(0,2.65)}{\pgfnodecenter{NodeW4}}}
   {\parbox{.65in}
   {\footnotesize Constructive mathematics}}{2pt}{2pt}
\pgfnodebox{NodeQ4}[stroke]
   {\pgfrelative{\pgfxy(3.3,0)}{\pgfnodecenter{NodeCD3}}}
   {{\footnotesize Con'ly}}{2pt}{2pt}  
\pgfnodebox{Therefore4}[virtual]
   {\pgfnodeborder{NodeQ4}{180}{8pt}}
   {\footnotesize So,}{2pt}{2pt} 
\pgfnodebox{NodeC4}[stroke]
   {\pgfrelative{\pgfxy(2,0)}{\pgfnodecenter{NodeQ4}}}
   {\parbox{.6in}{\footnotesize 
   $f(w)\not<m$,\\
   $f(w)\not>m$.}}{2pt}{2pt}
\pgfnodebox{Invisible5}[virtual]
   {\pgfnodeborder{NodeC4}{0}{6pt}}{}{0pt}{0pt} 
\pgfnodebox{NodeW5}[stroke]
   {\pgfrelative{\pgfxy(0,3.15)}{\pgfnodecenter{Invisible5}}}
   {\parbox{1.42in}
   {\footnotesize Trichotomy: If $x, y \in \mathbb{R}$ then $x < y$, $x = y$ or $x > y$.}}{2pt}{2pt}
\pgfnodebox{NodeB5}[stroke]
   {\pgfrelative{\pgfxy(0,2.65)}{\pgfnodecenter{NodeW5}}}
   {\parbox{.95in}{\footnotesize \centering Classical logic (specifically LEM)}}{2pt}{2pt}
\pgfnodebox{NodeQ5}[stroke]
   {\pgfrelative{\pgfxy(2.9,0)}{\pgfnodecenter{NodeC4}}}
   {{\footnotesize Class'ly}}{2pt}{2pt}  
\pgfnodebox{Therefore5}[virtual]
   {\pgfnodeborder{NodeQ5}{180}{8pt}}
   {\footnotesize So,}{2pt}{2pt} 
\pgfnodebox{NodeC5}[stroke]
   {\pgfrelative{\pgfxy(2.2,0)}{\pgfnodecenter{NodeQ5}}}
   {\parbox{.73in}{\footnotesize $f(w) = m$ for some $w$ s.t.  $u < w < v$.}}{2pt}{2pt}

\pgfputat{\pgfnodeborder{NodeC1}{90}{4pt}}{\pgfbox[center, base]{C$_{1}$}}
\pgfputat{\pgfnodeborder{NodeC3}{90}{4pt}}{\pgfbox[center, base]{C$_{2}$}}
\pgfputat{\pgfnodeborder{NodeW1}{0}{9pt}}{\pgfbox[center,center]{W$_{1}$}}
\pgfputat{\pgfnodeborder{NodeB1}{0}{8pt}}{\pgfbox[center,center]{B$_{1}$}}
\pgfputat{\pgfnodeborder{NodeW3}{0}{9pt}}{\pgfbox[center,center]{W$_{2}$}}
\pgfputat{\pgfnodeborder{NodeB3}{0}{8pt}}{\pgfbox[center,center]{B$_{2}$}}
\pgfputat{\pgfnodeborder{NodeW4}{0}{9pt}}{\pgfbox[center,center]{W$_{3}$}}
\pgfputat{\pgfnodeborder{NodeB4}{0}{8pt}}{\pgfbox[center,center]{B$_{3}$}}
\pgfputat{\pgfnodeborder{NodeQ1}{90}{4pt}}{\pgfbox[center,base]{Q$_{1}$}}
\pgfputat{\pgfnodeborder{NodeQ3}{90}{4pt}}{\pgfbox[center,base]{Q$_{2}$}}
\pgfputat{\pgfnodeborder{NodeQ4}{90}{4pt}}{\pgfbox[center,base]{Q$_{3}$}}
\pgfputat{\pgfnodeborder{NodeC4}{90}{4pt}}{\pgfbox[center,base]{C$_{3}$}}
\pgfputat{\pgfnodeborder{NodeD1}{90}{4pt}}{\pgfbox[center,base]{D$_{1}$}}
\pgfputat{\pgfnodeborder{NodeD3}{90}{4pt}}{\pgfbox[center,base]{D$_{3}$}}
\pgfputat{\pgfnodeborder{NodeW5}{0}{9pt}}{\pgfbox[center,center]{W$_{4}$}}
\pgfputat{\pgfnodeborder{NodeB5}{0}{8pt}}{\pgfbox[center,center]{B$_{4}$}}
\pgfputat{\pgfnodeborder{NodeQ5}{90}{4pt}}{\pgfbox[center,base]{Q$_{4}$}}
\pgfputat{\pgfnodeborder{NodeC5}{90}{4pt}}{\pgfbox[center, base]{C$_{4}$}}

\pgfsetlinewidth{1pt} 
\pgfsetendarrow{\pgfarrowsingle} 
\pgfnodesetsepstart{0pt}\pgfnodesetsepend{2pt} 
\pgfnodeconnline{NodeD1}{Therefore1} 
\pgfnodeconnline{NodeC1}{Therefore3} 
\pgfnodeconnline{NodeCD3}{Therefore4} 
\pgfnodeconnline{NodeQ1}{NodeC1} 
\pgfnodeconnline{NodeQ3}{NodeC3} 
\pgfnodeconnline{NodeQ4}{NodeC4} 
\pgfnodeconnline{NodeB1}{NodeW1} 
\pgfnodeconnline{NodeW1}{Invisible1} 
\pgfnodeconnline{NodeB3}{NodeW3} 
\pgfnodeconnline{NodeW3}{Invisible3} 
\pgfnodeconnline{NodeB4}{NodeW4} 
\pgfnodeconnline{NodeW4}{Invisible4} 
\pgfnodeconnline{NodeC4}{Therefore5} 
\pgfnodeconnline{NodeQ5}{NodeC5} 
\pgfnodeconnline{NodeB5}{NodeW5} 
\pgfnodeconnline{NodeW5}{Invisible5} 
\end{pgfmagnify} 
\end{pgfrotateby}
\end{pgfpicture} 
\end{center}
\caption{Classical proof of the Intermediate Value Theorem.}
\label{fig:ivtfull} \end{figure}

One source of examples of this phenomenon is the class of constructively invalid classical proofs.  Characteristically, most of the steps of these proofs \textit{are} constructively valid, making it important to identify the ones that are not.  Here a fine-grained application of Toulmin layouts can make the guilty steps explicit, unlike a single layout in which the qualifier for the whole proof would merely indicate that the result is classically, but not constructively, valid.  For instance, Figure~\ref{fig:root2} lays out a familiar classical proof that there are irrational numbers $\alpha$ and $\beta$ such that $\alpha ^{\beta}$ is rational.  (The rebuttal components have been omitted for simplicity.)  This decomposition of the proof into its separate steps clearly exhibits the dependencies of each step: whereas the second step relies on the constructively acceptable inference rule of constructive dilemma (CD), the first step employs the non-constructive law of excluded middle (LEM).  We can see immediately that the proof could be transformed into a constructive proof if a constructive derivation of $C_{1}$ from $D_{1}$ could be found to replace the first step.  This can be done, since there is a constructive proof that $\surd 2^{\surd 2}$ is irrational: substituting this statement for $W_{1}$, and its proof for $B_{1}$, yields a wholly constructive proof.\footnote{Albeit a rather cumbersome proof, since the constructive proof that $\surd 2^{\surd 2}$ is irrational relies on the deep Gelfond-Schneider Theorem.  There are much simpler constructive proofs that $\alpha ^{\beta}$ may be rational for irrational $\alpha$, $\beta$.}  
The Intermediate Value Theorem, which states that if $f$ is a continuous real-valued function, $u < v$ and $f(u) < m < f(v)$, then there must be a number $w$ such that  $u < w < v$ and $f(w) = m$, provides a more protracted example of a non-constructive classical proof.  The sketch of this proof laid out in Figure~\ref{fig:ivtfull}, has four steps, of which three are constructive.  However, the last step, that of Trichotomy, is irredeemably classical.  The constructivist has to choose between strengthening the hypotheses of the proof to rule out Brouwerian counterexamples to Trichotomy, and accepting a weaker result, that $f(w)$ and $m$ are arbitrarily close, rather than equal \cite[p.~140 f.]{George+}.  These examples suggest that, providing individual proof steps may be represented unambiguously, longer sequences of steps, whether represented separately or in combination, should also be unambiguous.

However, there are circumstances where ambiguity of layout would seem to be a genuine risk.  To begin with we shall look at another example from Euclid.  In 1879, Charles Dodgson, better known as Lewis Carroll, produced a diagram exhibiting the logical interdependency of the propositions of Book I \cite[frontispiece]{Carroll}.  A similarly motivated but distinct diagram for Book I (and diagrams for the other books, all inspired by the work of Ian Mueller \cite{Mueller}) is provided in the most recent scholarly edition of the \textit{Elements} \cite[p.\ 518]{Vitrac}.  Both diagrams omit some detail to maintain readability, and most of the differences between the two can be explained as contrasted editing choices.  However, there are some points at which they are explicitly at odds.  For example, for Mueller and Vitrac Proposition I.12 follows from I.8 and I.10, whereas for Carroll I.12 follows from I.9 alone.  Following Toulmin's practice of identifying the previously established propositions employed in a proof as its warrant, two separate layouts for the proof of I.12 could be constructed, with distinct warrants.  One possible explanation for the difference would be that Carroll was in error, or was working from a corrupt text: Johan Ludvig Heiberg's much improved edition of the \textit{Elements} only became available in Greek in 1883 and in English in 1908.  This text supports Mueller/Vitrac over Carroll \cite[Vol.~1, pp.\ 270 ff.]{Heath}.  If the difference can be resolved in this way it is obvious which layout is correct.

However, there are more persistent types of disagreement, one of which may provide an alternative account of this difference.  Any modeling process will make explicit ambiguities which had been ignored in the original context.  For example, propositional logic forces us to choose between \mbox{$A \wedge (B \vee C)$} and \mbox{$(A \wedge B) \vee C$} as formalizations of English sentences of the form ``$A$ and $B$ or $C$''.  We should be astonished to find an ambiguity of this kind in the \textit{statements} of a mathematical proof, since it would compromise the proof in a fashion any mathematician should be expected to spot.  Much less attention has been paid to the \textit{dialectic} of mathematical proofs.  In rendering it explicit through the formalism of the Toulmin layout we should not be surprised if we occasionally uncover a dialectical ambiguity: an argument which bears reconstruction in two distinct ways.  Where such ambiguities have gone uncorrected it is likely that they are sound arguments both ways.  This difference would be interesting, but benign.  If one or even both of the resolutions of a dialectical ambiguity is unsound, then the close reading required for the application of the Toulmin layout has performed a genuine service in demonstrating that the proof is at best poorly phrased and at worst fallacious. 

A more profound sort of disagreement may arise where there is no dispute over what the proof says, but there is a dispute over what it ought to say.   Some of the best known propositions of Euclid, such as I.5, the \textit{pons asinorum}, and I.47, Pythagoras's theorem, have many independent proofs, of which Euclid's is not necessarily the most rigorous, or the clearest.  Presumably, the different proofs will have different layouts, but there is no inconsistency in this. A more fundamental example from the \textit{Elements}, is the proof of Proposition I.4, that two triangles are equal if they have two sides and the enclosed angle equal, which makes use of the principle of superposition.  This principle, expressed in one of Euclid's axioms, Common Notion 4, as `things which coincide with one another are equal to one another', is used sparingly by Euclid, and has long been suspected as too empirical in character for geometrical proof \cite[Vol.~1, pp.\ 224 ff.]{Heath}.  Many modern axiomatizations of geometry take I.4 as an axiom, not a theorem \cite[Vol.~1, p.\ 249]{Heath}. \label{sec:I.4}

All of the sources of conflict considered so far reflect an ambiguity in the proof under analysis, either because the text of the proof is in dispute, or is ambiguously expressed, or because there are multiple distinct proofs of the same proposition.  All of these cases may give rise to a choice of different layouts, but in so doing they are faithfully reproducing an ambiguity in the source material.  This is a useful service, especially if a concealed ambiguity is brought to light.  Could a mathematical proof be represented by fundamentally distinct layouts in some other, harmful way?

\subsection{The Four Colour Theorem}

Alcolea, whose analysis of a critical mathematical argument I reproduced in Section \ref{sec:Reg}, also has a case study of a regular mathematical argument, Kenneth Appel and Wolfgang Haken's computer assisted proof of the four colour conjecture.  The conjecture states that four colours may be assigned to the regions of any planar map in such a way that no adjacent regions receive the same colour.  A partial proof published by Alfred Kempe in 1879 was taken as decisive for a decade until the discovery of counterexamples demonstrating its limitations.  The conjecture has proved resistant to straightforward methods of proof, and in 1976 was confirmed by a computer assisted proof, which remains controversial in some quarters since the full details are too protracted for human inspection.
Alcolea reconstructs the central argument of the proof as a derivation from the data $D_{1}$--$D_{3}$
\begin{quotation}
\noindent ($D_{1}$) Any planar map can be coloured with five colours.\\
($D_{2}$) There are some maps for which three colours are insufficient.\\
($D_{3}$) A computer has analysed every type of planar map and verified that each of them is 4-colorable.
\end{quotation}
of the claim $C$, that
\begin{quotation}
\noindent ($C$) Four colours suffice to colour any planar map.
\end{quotation}
by employment of the warrant $W$, 
\begin{quotation}
\noindent ($W$) The computer has been properly programmed and its hardware has no defects.
\end{quotation}
which has backing $B$
\begin{quotation}
\noindent ($B$) Technology and computer programming are sufficiently reliable. \cite[pp.~142f.]{Alcolea}
\end{quotation}
We can see at a glance that the warrant and backing of this proof are very different from those of Theaetetus's proof.  The dependence on apparently extra-mathematical methods is made explicit.  Alcolea draws the moral that the proof lacks mathematical rigour, and may even hide an unforeseen counterexample \cite[p.~143]{Alcolea}.

However, as I have argued elsewhere, this is not the only way of representing Appel and Haken's proof within Toulmin's layout \cite{Aberdein}.  We might also represent it as:
\begin{quotation}
Given that ($D$) the elements of the set $U$ are reducible, we can ($Q$) almost certainly claim that ($C$) four colours suffice to colour any planar map, since ($W$) $U$ is an unavoidable set (on account of ($B$) conventional mathematical techniques), unless ($R$) there has been an error in either (i) our mathematical reasoning, or (ii) the hardware or firmware of all the computers on which the algorithm establishing $D$ has been run.
\end{quotation}
This analysis requires some unpacking.  I have gone further into the details of the proof than Alcolea, in the hope of making the nature of its dependency on the computer precise.  The concepts of unavoidability and reducibility originate with Kempe's first, unsuccessful attempt to prove the conjecture.  An unavoidable set is a set of configurations---regions or clusters of neighbouring regions---such that every planar map must contain at least one member.  A configuration is reducible if it can be shown that all planar maps containing that configuration are four-colourable.  Kempe correctly demonstrated that the two-sided, three-sided, four-sided and five-sided regions together constitute an unavoidable set.  He also believed that he had shown all of these configurations to be reducible, which would have proved the conjecture.  Unfortunately, the five-sided region is not reducible.  Appel and Haken's proof is similar in nature to Kempe's, but much greater in scale: they derived an unavoidable set of reducible configurations with 1,482 members.  Although the set was found by a computer search, its unavoidability could be verified by hand.  However, confirmation of the reducibility of each of its members is too formidable a task for human proof checking.  Our confidence in the proof rests on the reliability of the methods programmed into the computer.

\subsection{Refutations}

The fundamental difference between these two layouts of the Four Colour Theorem lies in the rebuttal component, $R$.  Alcolea does not offer one, in accordance with Toulmin's injunction that `a mathematical inference \dots\ leaves no room for ``exceptions'' or ``rebuttals''; and indeed to raise the question of possible rebuttals would be to challenge the status of the entire argument' \cite[p.~254]{Toulmin+}.  On my reconstruction, the rebuttal forms the largest part of the layout.  What has happened here?

As the second part of this quotation from Toulmin makes clear, his account of rebuttal is not intended to permit challenges to the soundness of the argument, but rather to index possible exceptions allowed for by the choice of qualifier.  This distinction between the different ways in which an argument may be defeated has been developed at length in various places, often as a distinction between rebuttals and undercutters.\footnote{What follows is loosely derived from John Pollock's account.  He offers a general definition of a defeater as ``If $P$ is a reason for $S$ to believe $Q$, $R$ is a defeater for this reason if and only if $R$ is logically consistent with $P$ and $(P \wedge R)$ is not a reason for $S$ to believe $Q$'', and then defines rebutting and undercutting defeaters as:
\begin{quotation}
If $P$ is a \textit{prima facie} reason for $S$ to believe $Q$, $R$ is a \textit{rebutting} defeater for this reason if and only if $R$ is a defeater (for $P$ as a reason for $S$ to believe $Q$) and $R$ is a reason for $S$ to believe not-$Q$;

If $P$ is a \textit{prima facie} reason for $S$ to believe $Q$, $R$ is an \textit{undercutting} defeater for this reason if and only if $R$ is a defeater (for $P$ as a reason for $S$ to believe $Q$) and $R$ is a reason to deny that $P$ would not be true unless $Q$ were true \cite[pp.\ 38 f.]{Pollock}.
\end{quotation}}
Intuitively, an argument may be rebutted by offering independent reasons to disbelieve its conclusion, or undercut by challenging its soundness.  Most species of argument may be rebutted, undercut, or both.  In empirical science both defeaters occur, but as the result of very different work: an empirical study could be rebutted by an independent study, or undercut by a challenge to its methodology.  A paper which does the former might also do the latter, but need not.  Mathematics is different: proofs may be rebutted and undercut or just undercut, but not just rebutted.  Kempe's failed proof of the Four Colour Conjecture was undercut when counterexamples to it were discovered, but not rebutted since the counterexamples did not challenge the conjecture, which eventually turned out to be true.  Other failed proofs have been undercut and rebutted, when it has transpired that the `theorem' in question is not merely unproven but false.  But this is the only way a proof can be rebutted.  While we could accept that two empirical studies had conflicting results, despite both having sound methodologies, we could not accept that a mathematical conjecture could be independently proven to be both true and false: at least one of the `proofs' must fail.

Toulmin, however, does not admit undercutters into his layout.  Hence, the only sort of defeater he can accept is rebuttal without undercutting, which is precisely the sort of defeater which mathematical proofs cannot have.  In the light of his own definitions, Toulmin is right to say that mathematical proofs cannot be rebutted.  But should we accept these definitions?  Toulmin's motivation seems to be the thought that an undercutter cannot be part of a regular argument, but must instead reopen the issue `from a new, critical standpoint' \cite[p.~254]{Toulmin+}.  Is this what happens when a proof is undercut?  

Undercutters for mathematical proofs have been described as `inferential gaps' \cite[p.\ 51]{Fallis}.  In terms of my application of Toulmin's layout to multi-step proofs, an inferential gap would be a step in the proof in which the data and warrant do not support the claim.  The data represents the point the proof had reached when the (first) gap arises, the claim the point from which the remainder of the proof proceeds, and the warrant contains the mathematical results and definitions used to derive the claim from the data.  The failure of this inferential step might be attributable to the failure of the warrant.  Something like this happens when standards of rigour change, as in the criticism of Euclid's proposition I.4 discussed in Section \ref{sec:I.4}.  This would indeed open up a new critical argument.  However, what is far more typical of inferential gaps is that there is nothing wrong with the warrant---it just does not warrant the claim on the basis of the data.  It may be that the warrant needs to be supplemented, or that the claim is false, because the proof has been rebutted \textit{and} undercut, or has simply taken a wrong turning.  These are not grounds for a critical argument `challeng[ing] the credentials of current ideas', but evidence of human error in the conduct of a regular argument.  In the light of this observation, it seems more helpful to treat all defeaters alike, and admit them all into Toulmin's layout, generalizing his conception of `rebuttal'.\footnote{Several contemporary applications of Toulmin's work take this step.}

We can now understand the difference between the two reconstructions of the four colour theorem.  Alcolea sticks to Toulmin's narrow account of rebuttal, and produces a reconstruction with an emphatically non-mathematical warrant.  I offered instead a layout with a mathematical warrant, but a slightly more cautious qualifier than normal, admitting of two sources of rebuttal: human error and computer error.  I have argued elsewhere that in this case the likelihood of the former is orders of magnitude greater than the likelihood of the latter, and that this is why mathematicians are right to regard the proof with as much confidence as they do more orthodox proofs \cite{Aberdein}.

The underlying problem here is that of rational reconstruction, a familiar one in the history and philosophy of mathematics.  Mathematical textbooks unselfconsciously rewrite the proofs of past mathematicians, happily introducing anachronistic notation and standards of rigour.  Historians try to be more sensitive, but there is often a tension between doing justice to the context of discovery, by reproducing the twists and turns which led to the result, and doing justice to the context of justification by providing a sound proof \cite[p.\ 159]{Polya}.  This tension provides sufficient room for substantive disagreement about how best to reconstruct contentious proofs.  Toulmin's layout may bring these disagreements into sharper focus, but it can hardly be blamed for them.  As with the simpler sources of ambiguity discussed in Section \ref{sec:I.4}, the layout's capacity to represent proofs in different ways turns out not to be a handicap, but perhaps its most valuable feature.

\def\refname{\rm REFERENCES}
\setlength{\bibhang}{\parindent}

\end{document}